\documentclass[aps,preprint,amsmath,amssymb,amsfonts]{revtex4}
\usepackage{epsfig}
\usepackage{graphicx}
\usepackage{subfigure}
\usepackage{dcolumn}
\usepackage{bm}

\begin{document}

\title{Causal cells: spacetime polytopes with null hyperfaces}

\author{Yasha Neiman}
\email{yashula@gmail.com}
\affiliation{Institute for Gravitation \& the Cosmos and Physics Department, Penn State, University Park, PA 16802, USA}

\date{\today}

\begin{abstract}
We consider polyhedra and 4-polytopes in Minkowski spacetime - in particular, null polyhedra with zero volume, and 4-polytopes that have such polyhedra as their hyperfaces.
We present the basic properties of several classes of null-faced 4-polytopes: 4-simplices, ``tetrahedral diamonds'' and 4-parallelotopes. We propose a ``most regular'' representative of each class. The most-regular parallelotope is of particular interest: its edges, faces and hyperfaces are all congruent, and it features both null hyperplanes and null segments. A tiling of spacetime with copies of this polytope can be viewed alternatively as a lattice with null edges, such that each point is at the intersection of four lightrays in a tetrahedral pattern. We speculate on the relevance of this construct for discretizations of curved spacetime and for quantum gravity.
\end{abstract}
%
%
\maketitle

\section{Introduction}

Minkowski spacetime, the 4d pseudo-Euclidean space with metric $\eta_{\mu\nu} = \operatorname{diag}(-1,1,1,1)$, is central to the theory of Special Relativity. In the physical description of nature, it is arguably more important than the Euclidean plane or than 3d Euclidean space. Nevertheless, while planar and spatial geometry are richly developed, one is hard-pressed to find analogous results on the geometry of spacetime. Of course, the \emph{differential} geometry of spacetime is well-known to mathematicians and physicists alike. Yet there seems to be a shortage in ``school-level'' spacetime geometry, e.g. an analogue to the theory of Euclidean polygons and polyhedra. There is much intellectual satisfaction to be gained from this subject. One may also hope to arrive at physical insights, in particular with regard to discrete models of spacetime.

In this paper, we study the properties of some special 4-polytopes (the 4d version of polyhedra) in spacetime. The main qualitative difference between spacetime and Euclidean space is the existence of null (lightlike) directions. Thus, there exist line segments with vanishing length, plane elements with vanishing area, and hyperplane elements with vanishing volume. Of these, the 3d null hyperplane elements are especially interesting. In relativistic physics, null hypersurfaces play the role of causal boundaries between spacetime regions. They also function as characteristic surfaces for the differential equations of relativistic field theory. Important examples of null hypersurfaces include the lightcone of an event and the event horizon of a black hole. The prime example of a \emph{closed} null hypersurface is a causal diamond - the intersection of two lightcones originating from two timelike-separated points. Thus, our focus will be on 4-polytopes whose 3d hyperfaces are null polyhedra - linearly bounded regions of null hyperplanes. These 4-polytopes can be thought of as coarsely grained causal diamonds, with piecewise-flat boundaries. The 2d faces and 1d edges of such a polytope are necessarily spacelike (see section \ref{sec:polyhedra:general}). 

The paper is structured as follows. In section \ref{sec:hyperplanes}, we discuss the geometry of null hyperplanes. In section \ref{sec:polyhedra}, we discuss null polyhedra embedded in these hyperplanes. In section \ref{sec:polytopes}, we proceed to construct 4-polytopes out of these polyhedra. In section \ref{sec:polytopes:simplex}, we discuss the minimal null-faced polytope - the 4-simplex. In section \ref{sec:polytopes:diamond}, we discuss ``tetrahedral diamonds'' - shapes similar to causal diamonds, with spheres replaced by tetrahedra. In section \ref{sec:polytopes:parallel}, we discuss null-faced 4-parallelotopes. There we introduce the maximally regular 4-parallelotope. We refer to this shape as ``doubly null'', since in addition to having null hyperfaces, it contains a lightray segment as one of the diagonals of each hyperface. In section \ref{sec:tiling}, we discuss the straightforward tiling of spacetime with copies of the doubly-null parallelotope. This tiling has a dual description as a spacetime lattice threaded with lightrays in a tetrahedral pattern - in effect, a spacetime version of the diamond crystal lattice. As was brought to our attention during publication \footnote{I thank Ted Jacobson for informing me on this.}, both the tiling and its dual have been studied before, in \cite{Foster:2003jz} and \cite{Finkelstein:1996wu,Smith:1995kd} respectively, under the name ``hyperdiamond''.

The shapes we are considering are defined in flat Minkowski spacetime. However, flat polytopes and piecewise-flat tessellations also appear in discrete treatments of \emph{curved} spacetime, i.e. of General Relativity. On the classical level, this is the theory of Regge calculus \cite{Regge:1961px,Williams:1991cd}, where spacetime is constructed from flat simplices. In approaches to quantum gravity, simplicial triangulations appear prominently in causal dynamical triangulations \cite{Loll:2000my,Ambjorn:2010rx}. Somewhat less directly, flat polytopes emerge as an ingredient also in loop quantum gravity and spinfoam models \cite{Zakopane}. It is therefore of interest to consider curved tessellations constructed from our null-faced 4-polytopes. As we discuss in section \ref{sec:curved}, this is not feasible directly. Instead, we briefly speculate on adapting the lightray-threaded lattice of section \ref{sec:tiling} to a curved/gravitational setup.

We use $\mu,\nu$ as 4d spacetime indices. $a,b$ are 3d indices within a hyperplane. $\alpha,\beta$ are 2d indices within the space of lightrays in a null hyperplane. $i,j$ enumerate the faces of a polyhedron or the hyperfaces of a 4-polytope.

\section{Null hyperplanes} \label{sec:hyperplanes}

We wish to discuss null 3d polyhedra, or polyhedra with vanishing volume. These reside in null hyperplanes, such as the hyperplane $t=z$. Therefore, we will first briefly describe the geometry of these hyperplanes. The normal $\ell^\mu$ to the hyperplane ($\ell^\mu \sim (1,0,0,1)$ for the $t=z$ example) is a null vector, i.e. $\ell_\mu\ell^\mu = 0$. As a result, it is also tangent to the hyperplane. Its integral lines (lines of constant $x,y$ in the example) form null geodesics. The hyperplane is thus foliated into lightrays. All intervals within the hyperplane are spacelike, except the null intervals along the rays.

The hyperplane's intrinsic 3d metric $\gamma_{ab}$ is degenerate, with signature $(0,+,+)$. Since its determinant vanishes, all volumes within the hyperplane vanish. The intrinsic metric annihilates the null normal: $\gamma_{ab}\ell^b = 0$. It is useful to define a 2d quotient space $x^\alpha$ of the hyperplane's lightrays ($x^\alpha=(x,y)$ in the example), as opposed to the full 3d space $x^a$ of hyperplane points. The metric can then be viewed as a non-degenerate Euclidean metric $\gamma_{\alpha\beta} = \operatorname{diag}(1,1)$ on this quotient space. In other words, the metric measures distances \emph{between lightrays}, not between individual points. Arbitrary translations of shape elements along the lightrays do not affect the lengths and angles encoded by the metric. The determinant of the quotient metric $\gamma_{\alpha\beta}$ induces a multilinear \emph{area form}, as opposed to the volume form which exists in the non-null case. For a curved null hypersurface, the properties listed above hold locally, in the infinitesimal neighborhood of each point. 

Every lightray in flat spacetime defines a unique null hyperplane that contains it. Thus, when two null hyperplanes intersect, the intersection surface is always spacelike rather than null. This is also the case for curved null hypersurfaces.

In general, a degenerate metric $\gamma_{ab}$ does not induce a connection within the hypersurface, since it has no inverse. Thus, on a general null hypersurface, there is no notion of straight lines or of parallel transport. Instead, there is only a notion of parallel transport along the lightrays, which in particular supplies them with an affine structure. This limited structure is not derived from the intrinsic metric, but from the extrinsic curvature of the hypersurface. For a flat null hyperplane in Minkowski spacetime, more structure is available. There, we inherit from the ambient spacetime a full notion of parallel transport. As a result, we may speak of straight lines, compare the directions of line segments at different points, as well as compare the relative extents of lightray segments from different rays.

It is convenient to decompose the hyperplane coordinates as $x^a = (u,x^\alpha)$, where $x^\alpha$ labels the lightray, and $u$ is an affine coordinate along it. On a flat hyperplane, we can also match the scaling of $u$ between different rays. If in addition the $x^\alpha$ are 2d Cartesian coordinates, then $(u,x^\alpha)$ is the null analogue of a 3d Cartesian system. The hyperplane metric in such coordinates is $\gamma_{ab} = \operatorname{diag}(0,1,1)$. In the $t=z$ hyperplane, the choice $u = (t+z)/2$ and $x^\alpha = (x,y)$ answers these criteria. The symmetry group of the null hyperplane has a total of 7 generators, as follows:
\begin{itemize}
 \item Translations in the $x^\alpha$ plane (2 generators).
 \item Rotations in the $x^\alpha$ plane (1 generator).
 \item Translations in $u$ (1 generator).
 \item Boosts that rescale $u$ (1 generator): $u \rightarrow (1 + \epsilon)u$.
 \item Boosts that mix $u$ with the $x^\alpha$ (2 generators): $u \rightarrow u + \epsilon_\alpha x^\alpha$.
\end{itemize}
The symmetry group of a null hyperplane is larger by one generator than that of a non-null hyperplane. The reason is that the null coordinate $u$ does not have an overall scaling that must be preserved.

\section{Null polyhedra} \label{sec:polyhedra}

\subsection{General properties} \label{sec:polyhedra:general}

Within a null hyperplane, there exist both spacelike and null lines. Likewise, there exist spacelike and null planes. A null plane is a plane containing the null direction. The area form on such a plane vanishes. When constructing polyhedra in a null hyperplane, we will be using only spacelike edges and faces. This is because we are ultimately interested in 4-polytopes, so that the faces of the polyhedron should correspond to intersections with other null hyperplanes. As explained in the previous section, such intersections are necessarily spacelike. 

In 3d Euclidean space, each area element has a normal vector $\mathbf{n}$. When discussing polyhedra, it is convenient to define the norm of $\mathbf{n}$ to equal the area of the corresponding face. The orientation of the normals is chosen to be outgoing. Not every set of area normals $\left\{\mathbf{n}_i\right\}$ describes the faces of some polyhedron. For this to be true, the normals must sum up to zero:
\begin{align}
 \sum_i\mathbf{n}_i = 0 \ . \label{eq:vec_sum}
\end{align}
This can be understood as the requirement that the flux of any constant vector field through the polyhedron vanishes. In loop quantum gravity, the condition \eqref{eq:vec_sum} encodes the local $SO(3)$ rotation symmetry.

In a null hyperplane, a degenerate version of the above discussion holds. First, the normal direction to any area element within the hyperplane is the same. This is the null direction of the hyperplane's generating rays. One can still make the distinction between outgoing and ingoing face normals: the outgoing null direction is the one that takes us outside the polyhedron. Depending on the face in question, this may be either the future-pointing direction along the rays or the past-pointing one. The magnitude of the normal can no longer be tuned to the area of the face, since its metric length in any case vanishes. Instead, we can simply assign to each face its scalar area $s$, with a sign that indicates whether the outgoing normal is past-pointing or future-pointing. The analogue of eq. \eqref{eq:vec_sum} then reads:
\begin{align}
 \sum_i s_i = 0 \ . \label{eq:s_sum}
\end{align}
The condition \eqref{eq:s_sum} can be understood as the vanishing flux of an \emph{area current} $S^a$, which can be constructed from the 2d area density $\sqrt{\gamma}$. In analogy with the $SO(3)$ comment following \eqref{eq:vec_sum}, one can view eq. \eqref{eq:s_sum} as encoding an $SO(2)$ rotation symmetry in the 2d quotient space $x^\alpha$.

The distinction between ``past'' and ``future'' faces (those with past-pointing and future-pointing outgoing normals, respectively) is a useful one. The condition \eqref{eq:s_sum} can be expressed as an equality between the sum of ``past'' areas and the sum of ``future'' areas:
\begin{align}
 \sum s_+ = \sum s_- \ . \label{eq:s_diff}
\end{align}
In a convex polyhedron, the set of past faces is contiguous, as is the set of future faces. 

With regard to the combinatorics of vertices, edges and faces, a null polyhedron can be visualized just like an ordinary polyhedron in 3d Euclidean space. To correctly visualize its metric data, i.e. the lengths, areas and angles, a different approach is helpful. As explained in the previous section, the metric of a null hyperplane is essentially two-dimensional. For metric purposes, we can collapse each lightray to a point, leaving just the quotient plane $x^\alpha$. The edges and faces of the polyhedron become segments and polygons in this plane. Edges and faces that were ``above'' each other in the null direction now appear to intersect, and one must be careful to keep track of their identity. We will call this sort of picture the ``planar image'' of the polyhedron.

\subsection{Null tetrahedra} \label{sec:polyhedra:tetra}

\begin{figure}%
\centering%
\includegraphics[scale=1]{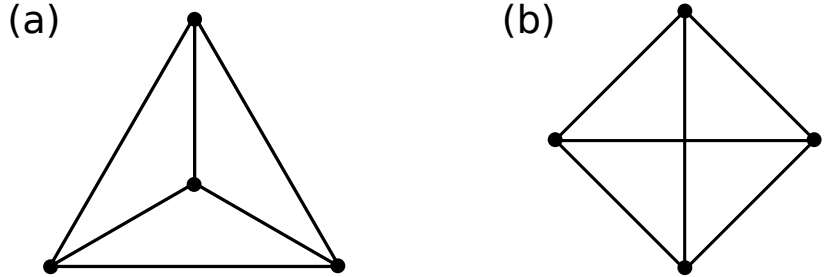} \\
\caption{The planar images of two null tetrahedra: (a) the maximally regular (1,3) tetrahedron; (b) the maximally regular (2,2) tetrahedron.}
\label{fig:tetra} 
\end{figure}

As in Euclidean space, the simplest null polyhedron is a tetrahedron. Up to reflections along the null axis, null tetrahedra come in two distinct types: (1,3) and (2,2). The pairs of numbers denote how many of the tetrahedron's four faces are past-pointing and future-pointing, respectively. Figure \ref{fig:tetra} shows the planar images of the two types of tetrahedron. This visualization makes the identity \eqref{eq:s_diff} obvious: the past-pointing and future-pointing faces literally occupy the same area. 

The shape of a null tetrahedron's planar image contains five degrees of freedom: $4\cdot 2 = 8$ vertex coordinates minus 3 isometries. It further turns out that the full shape of a null tetrahedron is completely determined by its planar image. Indeed, the $u$ coordinates of the four vertices can be tuned using the four independent symmetries $u \rightarrow \lambda u + \lambda_\alpha x^\alpha + u_0$. We note that these continuous symmetries cannot change causal relationships along a lightray, since the parameter $\lambda$ is positive. 

Since the null hyperplane contains a special direction, one cannot speak of regular null polyhedra in the full sense of the word. In figure \ref{fig:tetra}, we depicted the two ``maximally regular'' representatives of the two tetrahedron types. The tetrahedron in figure \ref{fig:tetra}(a) has as its faces an equilateral triangle with angles $(\pi/3,\pi/3,\pi/3)$ and three congruent isosceles triangles with angles $(2\pi/3,\pi/6,\pi/6)$. For the tetrahedron in figure \ref{fig:tetra}(b), the faces are four congruent isosceles triangles with angles $(\pi/2,\pi/4,\pi/4)$. As with regular tetrahedra in Euclidean space, each of the ``regular'' null tetrahedra is characterized by a single parameter: the overall metric scale of its planar image. This metric scale can be defined as the length of a given type of edge. 

\subsection{Null parallelepipeds} \label{sec:polyhedra:hexa}

After a tetrahedron, the next simplest polyhedron is a parallelepiped. The six faces of a null parallelepiped are spacelike parallelograms. There are three pairs of opposing faces, such that each pair is parallel and congruent. In a given pair of opposing faces, one is past-pointing, and the other future-pointing. The analogous statement in Euclidean space is that the outgoing normals of opposing faces point in opposite directions. Thus, the faces of a null parallelepiped are always arranged in a $(3,3)$ pattern: three past-pointing and three future-pointing. We will denote the intersection vertex of the three past faces as the ``initial'' vertex. Similarly, the three future faces intersect at the ``final'' vertex. Like a tetrahedron, a parallelepiped is fully defined by the positions of four vertices - for instance, the initial vertex and its three neighbors. Therefore, the shape of a null parallelepiped is also fully defined by its planar image, which again has five degrees of freedom.

\begin{figure}%
\centering%
\includegraphics[scale=0.75]{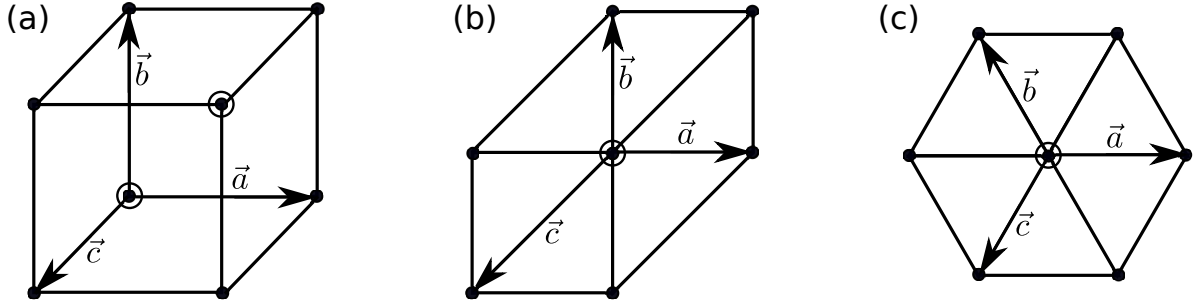} \\
\caption{The planar images of three null parallelepipeds: (a) a ``generic'' parallelepiped; (b) a doubly null parallelepiped; (c) the maximally regular doubly null parallelepiped. The initial and final vertices, which coincide in figures (b,c), are circled. The vectors $\mathbf{a},\mathbf{b},\mathbf{c}$ mark the three edges originating from the initial vertex.}
\label{fig:hexa} 
\end{figure}

The planar image of a null parallelepiped is shown in figure \ref{fig:hexa}(a). The initial and final vertices are situated inside a convex hexagon, formed by the other six vertices. The edges attached to the initial vertex divide the hexagon into three parallelograms, which are the three past faces. Similarly, the edges attached to the final vertex divide the same hexagon into the three future faces. This structure can be verified by construction, starting from the initial vertex and the three edges attached to it, and continuing by adding parallel and equal edges. A key point in the construction is that the initial vertex must be inside the triangle formed by its three neighbors. Otherwise, the planar images of two past faces would intersect. This would place two past faces directly ``above'' each other along the null axis, in violation of the parallelepiped's convexity. 

It is helpful to define three 2d vectors $\mathbf{a},\mathbf{b},\mathbf{c}$, giving the offsets from the initial vertex to its three neighbors in the planar image. Every edge in the planar image is described by one of these three vectors. The offset from the initial vertex to the final vertex is given by $\mathbf{a}+\mathbf{b}+\mathbf{c}$. Of particular interest is the case when $\mathbf{a}+\mathbf{b}+\mathbf{c} = 0$, depicted in figure \ref{fig:hexa}(b). In this case, the initial and final vertices coincide in the planar image. This means that they are situated on the same lightray in the null hyperplane. We will call such parallelepipeds ``doubly null'': in addition to being situated in a null hyperplane, they have a lightray segment as one of their diagonals. The shape of a doubly null parallelepiped is determined by three degrees of freedom. The other two are removed by the constraint $\mathbf{a}+\mathbf{b}+\mathbf{c} = 0$.

The most regular null parallelepiped is the doubly-null shape whose planar image is depicted in figure \ref{fig:hexa}(c). The six intermediate vertices form a regular hexagon. The initial and final vertices sit at the center of the hexagon, and the edges attached to them divide it into equilateral triangles. All twelve edges are of equal length, and all six faces are congruent rhombi with angles $(2\pi/3,\pi/3)$ (in figure \ref{fig:hexa}(c), each rhombus is given by the union of two equilateral triangles). This shape will play a special role in our discussion of 4-parallelotopes.

\section{Null-faced 4-polytopes} \label{sec:polytopes}

In this section, we construct spacetime 4-polytopes out of the null polyhedra discussed above. The geometric arena is expanded accordingly, from a single null hyperplane to the whole of spacetime. 

\subsection{Volume normals and the causal classification of faces} \label{sec:polytopes:class}

As with null polyhedra, the elements of spacetime 4-polytopes can be classified according to causality. In this paper, we are only considering null hyperfaces (and therefore, spacelike faces and edges). A null hyperplane bisects spacetime, with one half invariantly in the hyperplane's future, and the other in its past. Thus, we can classify every null hyperface as either ``past'' or ``future'', according to the direction from which it bounds the 4-polytope. 

In 4d Euclidean space, one can associate hyperfaces with outgoing ``volume normals'' - the equivalent of area normals in 3d. This notion is readily generalized for spacelike or timelike hyperfaces in spacetime. For null hyperfaces, the notion becomes more subtle. First, the normal to a null hyperface is tangent to it, so it can be neither outgoing nor ingoing. Second, the volume of a null hyperface vanishes. Nevertheless, a volume normal can still be usefully defined. Given three vectors $a^\mu,b^\mu,c^\mu$ in some null hyperplane, we define $\ell_\mu = \pm\epsilon_{\mu\nu\rho\sigma}a^\nu b^\rho c^\sigma / (3!)^2$ as the volume normal of the null tetrahedron spanned by $a^\mu,b^\mu,c^\mu$. The volume normal of an arbitrary null polyhedron can then be defined by triangulating it into tetrahedra and summing over them (alternatively, one could integrate over infinitesimal parallelepipeds). A completely analogous definition gives the customary volume normal in the non-null case. It remains to decide on the signs of the volume normals. This can be done as follows: choose the normal $\ell_\mu$ so that its scalar product with outgoing (necessarily non-normal) vectors is positive. This encodes the notion of an outgoing normal \emph{covector}, which is in fact more fundamental than that of a normal vector. With this choice, the volume normals satisfy a 4d version of the zero-sum identity \eqref{eq:vec_sum}:
\begin{align}
 \sum_i\ell^{(i)}_\mu = 0 \ . \label{eq:ell_sum}
\end{align}
With our mostly-plus convention for the metric signature, the above recipe associates past-pointing normal vectors $\ell^\mu$ to future hyperfaces, and future-pointing normals to past hyperfaces.   

While the magnitude of the null volume normals vanishes, the scalar products $\eta^{\mu\nu}\ell^{(i)}_\mu\ell^{(j)}_\nu$ of volume normals of \emph{different} hyperfaces are non-trivial. In a curved setting, these scalar products can be said to encode the inverse densitized spacetime metric $\sqrt{-g}g^{\mu\nu}$, in the same way that scalar products of edge vectors directly encode the metric $g_{\mu\nu}$. In addition, the scalar products between the volume normals of adjacent hyperfaces are a natural replacement for \emph{dihedral angles}. Note that a dihedral angle between null hyperplanes cannot be defined in the usual sense, as it corresponds to an infinite boost parameter. It is also worth noting that the scalar product of two null vectors $\ell^{(i)}_\mu,\ell^{(j)}_\mu$ coincides in magnitude with the area of the timelike parallelepiped spanned by them:
\begin{align}
 \left(\frac{\ell^{(i)}_\mu\ell^{(j)}_\nu - \ell^{(i)}_\nu\ell^{(j)}_\mu}{\sqrt{2}}\right)^2 = -\eta^{\mu\nu}\ell^{(i)}_\mu\ell^{(j)}_\nu \ . \label{eq:timelike_area}
\end{align}

\begin{figure}%
\centering%
\includegraphics[scale=1]{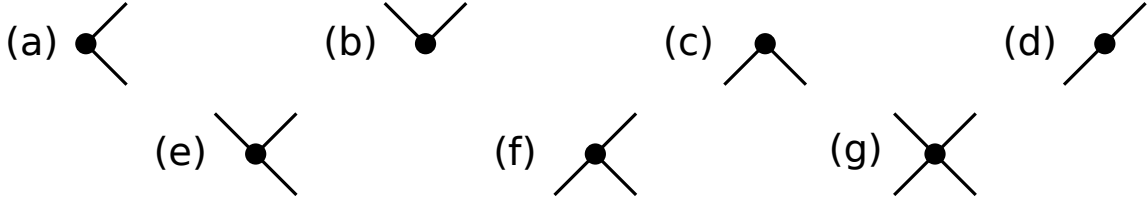} \\
\caption{Different choices of legs on the lightcross. Up and down stand for future and past, respectively. There is no distinction between right and left, apart from telling which legs belong to the same ray. Figures (a,b,c) depict past-future, past-past and future-future 2d faces of a 4-polytope, respectively. Figure (d) depicts a null ray/hyperplane continuing straight through the plane element. Figures (e,f,g) are relevant for tessellations of spacetime with null hyperfaces.}
\label{fig:lightcross} 
\end{figure}

Each 2d face of a 4-polytope is an intersection of two hyperfaces. Thus, the 2d faces also have a causal classification, according to the types of their associated hyperfaces: ``past-past'', ``past-future'' and ``future-future''. There is another fruitful way to view the different types of faces. A spacelike plane element in spacetime is orthogonal to a timelike plane, which contains exactly two lightrays. Thus, there is a ``lightcross'' of two lightrays orthogonal to the plane element. The lightcross has four ``legs'', each a half-lightray. Now, these half-rays determine the directions in which the plane element can be extended into a null hyperplane - the unique hyperplane containing the chosen ray. Therefore, when speaking of two null hyperfaces intersecting at a face, we must specify which two legs of the lightcross they occupy. This results in the same classification into ``past-past'', ``past-future'' and ``future-future'' faces. See figure \ref{fig:lightcross}(a,b,c).

The edges and vertices can be similarly classified according to the types of hyperfaces that contain them. Thus, if an edge or a vertex is at the intersection of three past and one future hyperface, we can classify it as past-past-past-future. This classification is less universal than that of faces and hyperfaces, since the number of hyperfaces intersecting at an edge or a vertex may vary.

\subsection{Null-faced 4-simplices} \label{sec:polytopes:simplex}

The simplest class of polytopes are simplices. Accordingly, as our first example of 4-polytopes with null hyperfaces, we discuss 4-simplices whose hyperfaces have zero volume. A 4-simplex has five tetrahedral hyperfaces, which in our case will be null tetrahedra, described in section \ref{sec:polyhedra:tetra}. A simple way to construct a null-faced 4-simplex is to draw five null hyperplanes in spacetime. The 4-simplex is then defined by the convex hull of their intersections. The five tetrahedra all intersect each other, with a total of 10 triangular faces. The number of edges is also 10, and the number of vertices is 5. 

\begin{figure}%
\centering%
\includegraphics[scale=1]{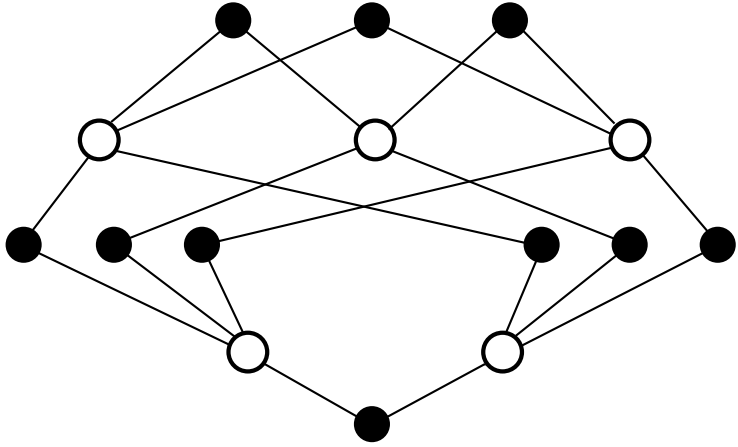} \\
\caption{The causal structure of a null-faced 4-simplex (up to time reversal). The axes are as in figure \ref{fig:lightcross}. The empty circles stand for tetrahedral hyperfaces, while the full circles stand for triangular faces. A link between a face and a hyperface indicates that one is contained in the other.}
\label{fig:simplex_causal} 
\end{figure}

It follows from convexity that out of the five tetrahedra, two must act as past hyperfaces and three as future hyperfaces, or vice versa. The two past tetrahedra are of (1,3) type, while the three future ones are of (2,2) type, with the obvious replacements for the time-reversed case. Of the ten triangular faces, one is past-past type, six are past-future, and three are future-future. Of the ten edges, three are past-past-future, six are past-future-future, and one is future-future-future. Of the five vertices, three are past-past-future-future, and two are past-future-future-future. Every pair of vertices is spacelike-separated, since they are connected by an edge. The causal structure of the faces and hyperfaces is depicted in figure \ref{fig:simplex_causal}.

One consequence of this causal structure is that spacetime \emph{cannot} be triangulated by null-faced 4-simplices (though, of course, it \emph{can} be triangulated by more general 4-simplices). The reason is that in any tessellation of spacetime into convex shapes, the 2d faces must connect to hyperfaces according to one of the structures (e,f,g) in figure \ref{fig:lightcross}. In all of these, the number of past-past plus the number of future-future ``dihedral angles'' equals the number of past-future ones. On the other hand, we've seen that in any null-faced 4-simplex, the number of past-future angles is larger: 6, as opposed to $1+3=4$ past-past and future-future angles.

The shape of a null-faced 4-simplex is determined by five degrees of freedom. These can be represented in a number of useful ways:
\begin{itemize}
 \item $5\cdot 4 = 20$ coordinates for the five vertices, minus 5 zero-volume constraints on the hyperfaces, minus 10 isometries.
 \item 10 edge lengths, minus 5 zero-volume constraints.
 \item 10 face areas, minus 5 zero-sum constraints of the form \eqref{eq:s_diff}. These are equivalent to the zero-volume constraints in the previous descriptions.
 \item $5\cdot 2 = 10$ degrees of freedom for the directions of the five null hyperplanes, plus 5 degrees of freedom for their offsets, minus 10 isometries.
 \item $5\cdot 3 = 15$ degrees of freedom for the five null volume normals $\ell^{(i)}_\mu$, minus 4 components of the zero-sum constraint \eqref{eq:ell_sum}, minus 6 rotations.
 \item $5\cdot 2 = 10$ degrees of freedom for the directions of the five null hyperfaces, minus 6 rotations, plus 1 scale parameter, such as the spacetime volume $\Omega$ of the 4-simplex.
 \item $5$ degrees of freedom for the shape of a single tetrahedral null hyperface. Thanks to the lightcross structure, the position and orientation of its four faces uniquely determine the other four null hyperplanes.
\end{itemize}

The scalar products $\eta^{\mu\nu}\ell_\mu^{(i)}\ell_\nu^{(j)}$ of the null volume normals are directly related to the spacetime volume of the 4-simplex and to the areas of the 2d faces. To express the spacetime volume, we must choose a set of four volume normals $\ell_\mu^{(i)}$. The time-orientation of the normals should be correlated with the past/future status of their hyperfaces, as prescribed before eq. \eqref{eq:ell_sum}. Next, we construct a symmetric $4\times 4$ matrix $L^{ij} = (3!)^2\eta^{\mu\nu}\ell_\mu^{(i)}\ell_\nu^{(j)}$ of their scalar products. The diagonal elements of $L^{ij}$ are zero. Elements corresponding to past-future pairs $ij$ are positive, while those for past-past and future-future pairs are negative. The spacetime volume can then be found as:
\begin{align}
 \Omega = \frac{1}{4!}\left|\det L\right|^{1/6} \ . \label{eq:Omega_L_simplex}
\end{align}
The area of the face at the intersection of the $i$'th and $j$'th hyperplanes can be found as:
\begin{align}
 s^{ij} = \frac{(3!)^2\left|\eta^{\mu\nu}\ell_\mu^{(i)}\ell_\nu^{(j)}\right|}{2!4!\Omega} = \frac{\left|L^{ij}\right|}{2\left|\det L\right|^{1/6}} \ ,
 \label{eq:s_simplex}
\end{align}
where $L^{ij}$ is constructed from a set of four volume normals that includes the desired pair $ij$. While eq. \eqref{eq:Omega_L_simplex} has an immediate analogue for general 4-simplices, eq. \eqref{eq:s_simplex} is unique to the null case. It arises from the equivalence \eqref{eq:timelike_area} between scalar products and timelike areas. One can combine the properties \eqref{eq:Omega_L_simplex}-\eqref{eq:s_simplex} and define the 4-volume directly in terms of triangle areas:
\begin{align}
 \Omega = \frac{1}{4!}\left|\det(2S)\right|^{1/2} = \frac{1}{6}\left|\det S\right|^{1/2} \ . \label{eq:Omega_S_simplex}
\end{align}
Here, $S^{ij}$ is a symmetric $4\times 4$ matrix of face areas, with zeros on the diagonal. The off-diagonal elements are $S^{ij} = \pm s^{ij}$, with plus for past-future faces, and minus for past-past and future-future ones.

\begin{figure}%
\centering%
\includegraphics[scale=1]{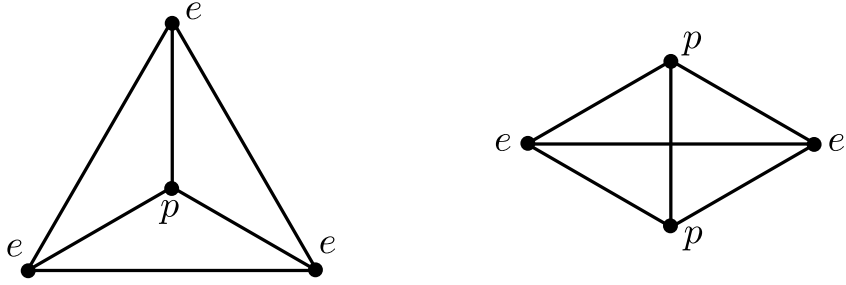} \\
\caption{The two types of tetrahedral hyperfaces on the maximally symmetric 4-simplex. The $e$/$p$ labels correspond to the two types of vertices in eq. \eqref{eq:simplex_vertices}.}
\label{fig:tetra_simplex} 
\end{figure}

As should be clear from the causal structure in figure \ref{fig:simplex_causal}, a fully regular null-faced 4-simplex does not exist: the past and future hyperfaces are necessarily different. Also, the two past tetrahedra (in the time-orientation of figure \ref{fig:simplex_causal}) are of (1,3) type, i.e. they have one face whose area equals the sum of the three others. Thus, we must have at least two different face areas, and therefore at least two different edge lengths. There exists a 4-simplex with \emph{precisely} two different edge lengths and two different face areas. We will call it the maximally regular null-faced 4-simplex. This shape has a single free parameter, determining its overall scale. Perhaps the most elegant way to define this 4-simplex is by specifying the directions of its five null volume normals. These directions can be viewed as points on a Riemann sphere \cite{Penrose:1985jw}. We place the two normals to the past hyperfaces at the north and south poles of the sphere. The normals to the three future hyperfaces are then placed at the vertices of an equilateral triangle on the equator. In appropriate axes $(t,x,y,z)$, the resulting 4-simplex has the following vertices:
\begin{align}
 & e_1: (0,-\sqrt{3}a,-a,0); && e_2: (0,\sqrt{3}a,-a,0); && e_3: (0,0,2a,0); \nonumber \\
 & p_1: (a,0,0,-a); && p_2: (a,0,0,a) \ . \label{eq:simplex_vertices}
\end{align}
Here, $e$ stands for ``equator'' vertices, while $p$ stands for ``pole'' vertices. The three $(ee)$-type edges have length $2\sqrt{3}a$, while the $(pp)$ edge and the six $(ep)$ edges have length $2a$. The $(eee)$ triangle and the three smaller $(epp)$ triangles are equilateral, with angles $(\pi/3,\pi/3,\pi/3)$. The six $(eep)$ triangles are isosceles, with angles $(2\pi/3,\pi/6,\pi/6)$. The $(eee)$ triangle has area $3\sqrt{3}a^2$, while the $(epp)$ and $(eep)$ triangles have area $\sqrt{3}a^2$. The two $(eeep)$ hyperfaces are maximally regular (1,3)-type tetrahedra. The three $(eepp)$ tetrahedra are of (2,2) type. The planar images of the two types of hyperfaces are depicted in figure \ref{fig:tetra_simplex}. The spacetime volume of the 4-simplex is $\Omega = \sqrt{3}a^4/2$.

\subsection{Tetrahedral diamonds} \label{sec:polytopes:diamond}

The next polytope we will discuss is perhaps the most similar to a piecewise-flat causal diamond. Instead of a sphere at the intersection of two lightcones, this shape is based on a spacelike tetrahedron. Thus, we begin with an arbitrary spacelike tetrahedron, situated without loss of generality in the $t=0$ hyperplane. We will call this the ``base tetrahedron''. For each of the base tetrahedron's four faces, we draw the ``lightcross'' of two null hyperplanes orthogonal to it. Our ``tetrahedral diamond'' is then defined by the convex hull of the intersections of these null hyperplanes. 

The four past-ingoing hyperplanes all intersect at a single point, as do the four future-ingoing ones. These ``initial'' and ``final'' vertices are situated directly to the past and future (along the $t$ axis) from the center of the base tetrahedron's inscribed sphere. To see this, note that the four past hyperfaces are tangent to the future lightcone of the initial vertex, and likewise for the future hyperfaces. The intersection of the two lightcones is therefore the inscribed sphere of the base tetrahedron. Placing the center of the inscribed sphere at the origin and denoting its radius by $r$, the coordinates of the initial and final vertices are $(\pm r,0,0,0)$. 

Overall, the tetrahedral diamond has six vertices: the four vertices of the base tetrahedron, plus the initial vertex and the final vertex. The number of edges is 14: the six edges of the base tetrahedron, plus four connecting its four vertices to the initial vertex, plus four connecting to the final vertex. The 2d faces are 16 triangles: the four faces of the base tetrahedron, plus six triangles connecting the edges of the base tetrahedron to the initial vertex, plus six triangles connected to the final vertex. The hyperfaces are eight null tetrahedra: four tetrahedra connecting the faces of the base tetrahedron to the initial vertex, plus four connected to the final vertex. 

The tetrahedral diamond can be viewed as the union of two congruent 4-simplices, one ``initial'' and one ``final'', glued to each other along the base tetrahedron. The five hyperfaces of each 4-simplex consist of the spacelike base tetrahedron and the four past (future) null tetrahedra connecting its faces to the initial (final) vertex.

Out of the eight hyperfaces of the tetrahedral diamond, four are past-type (the ones connected to the initial vertex), and four are future-type (the ones connected to the final vertex). The past and future tetrahedra are of (3,1) and (1,3) types, respectively. Of the 16 triangular faces, six are past-past, four are past-future, and six are future-future. Of the 14 edges, four are past-past-past, six are past-past-future-future, and four are future-future-future. Of the 6 vertices, the initial vertex is past-past-past-past, the final one is future-future-future-future, and the four base vertices are past-past-past-future-future-future. The initial and final vertices are timelike-separated. Every other pair of vertices is spacelike-separated, since they share a 2d face.

\begin{figure}%
\centering%
\includegraphics[scale=1]{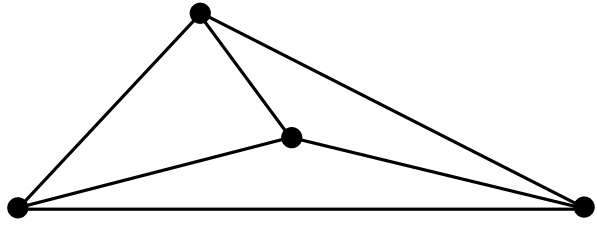} \\
\caption{A triangular face of the spacelike base tetrahedron. In the interior, we marked the tangency point between the face and the tetrahedron's inscribed sphere. The same figure can be interpreted as the planar image of the past or future null tetrahedron connected to the chosen face.}
\label{fig:tetra_face} 
\end{figure}

In figure \ref{fig:tetra_face}, we depict one triangular face of the base tetrahedron. The interior point connected to the vertices represents the point of tangency between the face and the base tetrahedron's inscribed sphere. The segments connecting this point to the initial and final vertices are null. As a result, figure \ref{fig:tetra_face} can be immediately interpreted as the planar image of the (congruent) past and future null tetrahedra connected to the chosen face. This observation reduces the task of calculating the edge lengths and face areas of the tetrahedral diamond to the 3d geometry of the spacelike tetrahedron.

The 4-volume of a tetrahedral diamond can be found as twice the volume of a 4-simplex, with the spacelike tetrahedron as its base and the inscribed radius $r$ as its height. The result is:
\begin{align}
 \Omega = 2\cdot\frac{1}{4}rV = \frac{1}{2}rV \ , \label{eq:Omega_diamond}
\end{align}
where $V$ is the base tetrahedron's volume. As with null 4-simplices, it's impossible to tessellate spacetime with tetrahedral diamonds. The reason is again a mismatch between the numbers of 2d faces of different causal types. In this case, the mismatch is in the opposite direction: the number of past-future faces (four) is smaller than the number of past-past plus future-future faces ($6+6 = 12$).

The ``maximally regular'' tetrahedral diamond is of course constructed from a regular base tetrahedron. The resulting shape has a single free parameter, determining its overall scale. In appropriate axes $(t,x,y,z)$, the six vertices read:
\begin{align}
 & i: (-a,0,0,0); \nonumber \\
 & b_1: (0,-\sqrt{6}a,-\sqrt{2}a,-a); && b_2: (0,\sqrt{6}a,-\sqrt{2}a,-a); && b_3: (0,0,2\sqrt{2}a,-a); && b_4: (0,0,0,3a); \nonumber \\
 & f: (a,0,0,0) \ .
\end{align}
Here, $b$ stands for ``base'', $i$ for ``initial'', and $f$ for ``final''. The inscribed radius of the base tetrahedron is $r = a$. The lengths of the base tetrahedron's six edges are $2\sqrt{6}a$. The lengths of the eight edges connecting it to the initial and final vertices are $2\sqrt{2}a$. The areas of the base tetrahedron's four faces are $6\sqrt{3}a^2$. The areas of the twelve faces connecting it to the initial and final vertices are $2\sqrt{3}a^2$. The past/future null hyperfaces are all congruent maximally-regular (1,3)/(3,1) tetrahedra, with planar images as in figure \ref{fig:tetra}(a) or \ref{fig:tetra_simplex}(a). The 4-volume of the tetrahedral diamond is $\Omega = 4\sqrt{3}a^4$.

\subsection{Null-faced parallelotopes} \label{sec:polytopes:parallel}

The last class of polytopes that we discuss are spacetime parallelotopes with null hyperfaces. As with any parallelotope, copies of a null-faced parallelotope generate a tiling of Minkowski spacetime. This sets it apart from the shapes we considered previously. 

One way to construct a null-faced parallelotope is to choose four linearly independent pairs of parallel null hyperplanes. The parallelotope is the convex shape demarcated by the hyperplanes' intersections. A 4d parallelotope has 16 vertices. It has 32 edges, in four groups of eight, where each group consists of parallel edges of equal length. The 24 faces are spacelike parallelograms, divided into six groups of four, with the faces in each group parallel and congruent. The 8 hyperfaces are null parallelepipeds, described in section \ref{sec:polyhedra:hexa}. These are divided into four opposing pairs, with the hyperfaces in each pair parallel and congruent. 

In each opposing pair of hyperfaces, one is past-type, while the other is future-type. Thus, there are four past hyperfaces, and four future ones. The four past hyperfaces intersect at an ``initial'' vertex, while the four future hyperfaces intersect at a ``final'' vertex. The 2d faces are divided into six past-past (connecting every pair of past hyperfaces), six future-future (connecting every pair of future hyperfaces), and twelve past-future (connecting every non-opposing pair of past and future hyperfaces). The 32 edges are divided into four past-past-past (connected to the initial vertex), four future-future-future (connected to the final vertex), twelve past-past-future and twelve past-future-future edges. The 16 vertices are divided into one past-past-past-past (the initial vertex, which is also the initial vertex of the four past hyperfaces), one future-future-future-future (the final vertex, also the final vertex of the four future hyperfaces), four past-past-past-future (the initial vertices of the four future hyperfaces), four past-future-future-future (the final vertices of the four past hyperfaces), and six past-past-future-future.

As usual, the properties of parallelotopes closely mirror the properties of simplices. Like a 4-simplex, a 4d parallelotope is fully defined by the positions of five vertices, such as the initial vertex and its four neighbors. However, there is a difference in the null-faced case: while a null-faced 4-simplex comes with five independent zero-volume constraints, a null-faced parallelotope has only four. As a result, the shape of a null-faced parallelotope is determined by six degrees of freedom - one more than for a null-faced 4-simplex. These degrees of freedom can be represented in the following ways:
\begin{itemize}
 \item $5\cdot 4 = 20$ coordinates for the initial vertex and its four neighbors, minus 4 zero-volume constraints on the independent hyperfaces, minus 10 isometries.
 \item 4 lengths plus 6 angles among the independent edges, minus 4 zero-volume constraints.
 \item 6 independent face areas.
 \item $4\cdot 2 = 8$ degrees of freedom for the directions of the past null hyperfaces, plus 4 degrees of freedom for the offsets between opposing pairs, minus 6 rotations.
 \item $4\cdot 3 = 12$ degrees of freedom for the null volume normals $\ell^{(i)}_\mu$ to the past hyperfaces, minus 6 rotations.
 \item $5$ degrees of freedom for the shape of a single null parallelepiped, plus 1 for the length of the edges that connect it to the opposing hyperface. Thanks to the lightcross structure, the position and orientation of the parallelepiped's faces uniquely determine the other six null hyperplanes.
\end{itemize}

As with the 4-simplex, the scalar products $\eta^{\mu\nu}\ell_\mu^{(i)}\ell_\nu^{(j)}$ of four independent null volume normals (e.g. the normals to the past hyperfaces) can be related to the 4-volume and to the face areas. Similar formulas hold, but without the factorials. Again, we should correlate the time-orientation of the normals with the past/future status of their hyperfaces. We then construct the symmetric $4\times 4$ matrix $L^{ij} = \eta^{\mu\nu}\ell_\mu^{(i)}\ell_\nu^{(j)}$ of scalar products, with zeroes on the diagonal. Elements corresponding to past-future pairs $ij$ are positive, while those for past-past and future-future pairs are negative. The spacetime volume is then given by:
\begin{align}
 \Omega = \left|\det L\right|^{1/6} \ . \label{eq:Omega_L_para}
\end{align}
The area of the face at the intersection of the $i$'th and $j$'th hyperplanes can be found as:
\begin{align}
 s^{ij} = \frac{\left|\eta^{\mu\nu}\ell_\mu^{(i)}\ell_\nu^{(j)}\right|}{\Omega} = \frac{\left|L^{ij}\right|}{\left|\det L\right|^{1/6}} \ .
 \label{eq:s_para}
\end{align}
Combining eqs. \eqref{eq:Omega_L_para}-\eqref{eq:s_para}, we can define the 4-volume directly in terms of parallelogram areas:
\begin{align}
 \Omega = \left|\det S\right|^{1/2} \ . \label{eq:Omega_S_para}
\end{align}
Here, $S^{ij}$ is a symmetric $4\times 4$ matrix of face areas, with zeroes on the diagonal. The off-diagonal elements are $S^{ij} = \pm s^{ij}$, with plus for past-future faces, and minus for past-past and future-future ones. If we concentrate on the four past hyperfaces (or the four future ones) as our independent set, all the sign factors in the above can be disregarded.

Of particular interest are null-faced parallelotopes whose hyperfaces are all \emph{doubly-null} parallelepipeds, as defined in section \ref{sec:polyhedra:hexa} and depicted in figure \ref{fig:hexa}(b,c). It turns out that up to an overall scale, there is exactly one shape of this kind. Indeed, recall that a doubly-null parallelepiped obeys the constraint $\mathbf{a}+\mathbf{b}+\mathbf{c} = 0$ on the planar images of its edges. This means that the lengths of two independent edges and the angle between them determine the length and direction of the third edge. Now, out of the four past (say) hyperfaces of the parallelotope, every pair intersects at a 2d parallelogram, thus sharing two edges and the angle between them. This results in a web of equality constraints among the different edges, whose only solution is to set all the edge lengths equal. Thus, the hyperfaces of a doubly-null 4-parallelotope are all congruent maximally-regular doubly-null parallelepipeds, depicted in figure \ref{fig:hexa}(c).

The doubly-null 4-parallelotope may be called the most regular parallelotope in spacetime. All its edges, all its faces and all its hyperfaces are congruent. The 2d faces are rhombi with angles $(2\pi/3,\pi/3)$. The hyperfaces are maximally regular doubly-null parallelepipeds, i.e. their initial and final vertices are connected by lightray segments. To write down the vertices of this parallelotope, we choose an arbitrary ordering of the four opposing pairs of hyperfaces. The vertices can then be denoted by a 4-tuple of $p$'s and $f$'s, to indicate whether the vertex belongs to the past or future hyperface within each pair. In appropriate axes $(t,x,y,z)$, the vertex coordinates then read:
\begin{align}
 & pppp: (-a,0,0,0); \nonumber \\
 & fppp: \left(-\frac{a}{2},+\frac{\sqrt{3}a}{2},+\frac{\sqrt{3}a}{2},+\frac{\sqrt{3}a}{2}\right); \quad
 pfpp: \left(-\frac{a}{2},+\frac{\sqrt{3}a}{2},-\frac{\sqrt{3}a}{2},-\frac{\sqrt{3}a}{2}\right); \quad \nonumber \\
 & ppfp: \left(-\frac{a}{2},-\frac{\sqrt{3}a}{2},+\frac{\sqrt{3}a}{2},-\frac{\sqrt{3}a}{2}\right); \quad  
 pppf: \left(-\frac{a}{2},-\frac{\sqrt{3}a}{2},-\frac{\sqrt{3}a}{2},+\frac{\sqrt{3}a}{2}\right); \quad \nonumber \\
 & ffpp: \left(0,+\sqrt{3}a,0,0\right);\quad fpfp: \left(0,0,+\sqrt{3}a,0\right);\quad fppf: \left(0,0,0,+\sqrt{3}a\right); \nonumber \\
 & ppff: \left(0,-\sqrt{3}a,0,0\right);\quad pfpf: \left(0,0,-\sqrt{3}a,0\right);\quad pffp: \left(0,0,0,-\sqrt{3}a\right); \nonumber \\
 & pfff: \left(+\frac{a}{2},-\frac{\sqrt{3}a}{2},-\frac{\sqrt{3}a}{2},-\frac{\sqrt{3}a}{2}\right); \quad
 fpff: \left(+\frac{a}{2},-\frac{\sqrt{3}a}{2},+\frac{\sqrt{3}a}{2},+\frac{\sqrt{3}a}{2}\right); \quad \nonumber \\
 & ffpf: \left(+\frac{a}{2},+\frac{\sqrt{3}a}{2},-\frac{\sqrt{3}a}{2},+\frac{\sqrt{3}a}{2}\right); \quad  
 fffp: \left(+\frac{a}{2},+\frac{\sqrt{3}a}{2},+\frac{\sqrt{3}a}{2},-\frac{\sqrt{3}a}{2}\right); \quad \nonumber \\
 & ffff: (+a,0,0,0) \ . \label{eq:para_vertices}
\end{align}
The edge length is $\sqrt{2}a$. The area of each face is $\sqrt{3}a^2$. The spacetime volume is $3\sqrt{3}a^4$. The initial and final vertices of the parallelotope ($pppp$ and $ffff$) are timelike-separated. The initial and final vertices of each 3d face ($pppp$ and $pfff,fpff,ffpf,fffp$ for the past faces, $ffff$ and $fppp,pfpp,ppfp,pppf$ for the future faces) are null-separated. Every other pair of vertices is spacelike-separated. Each causal ``level'' of vertices (defined by the number of $p$'s vs. the number of $f$'s) sits in a constant-time hyperplane, with equal spacing $\Delta t = a/2$ among the different levels. The four past-past-past-future vertices and the four past-future-future-future vertices each form a regular tetrahedron, with edge length $\sqrt{6}a$. The six past-past-future-future vertices form a regular octahedron, also with edge length $\sqrt{6}a$.

\section{Tiling spacetime with the doubly-null parallelotope} \label{sec:tiling}

Like any parallelotope, the doubly-null parallelotope of eq. \eqref{eq:para_vertices} generates a tiling of spacetime. In this tiling, each face is attached to a lightcross of null hyperfaces, as depicted in figure \ref{fig:lightcross}(g). The tiling is remarkable for its regularity. In particular, unlike the naive cubic tiling along the $(t,x,y,z)$ axes, its edges, faces and hyperfaces are all congruent, without distinct spacelike and timelike elements. We note, however, that the tiling does have a preferred reference frame - the one in which the initial and final vertices of the 4-parallelotope are directly above each other along the time axis.

Another feature of this tiling is due to the null diagonals of the hyperfaces. Each node in the tiling acts as the initial vertex of one 4-parallelotope, and as the final vertex of another. It is therefore connected by null segments to eight other nodes - four to the future, along the past faces of the first parallelotope, and four to the past, along the future faces of the second one. These eight null segments sit on four lightrays that pass through the original node, one segment along each half-ray. As should be clear from regularity, the directions of these four rays form a tetrahedral pattern (this is in fact a statement about the cross-ratio of the four null directions \cite{Penrose:1985jw}). 

We can now strip away the edges, faces and hyperfaces of the original tiling, leaving only the lattice of nodes and the lightrays that connect them. This dual picture describes a ``cubic'' (in fact, parallelotopic) lattice threaded by lightrays. There are four different directions of rays, forming a tetrahedral pattern. Four rays pass through each lattice node. Along each ray, an infinite number of affinely-spaced lattice nodes are threaded. Due to the regular lattice structure and the tetrahedral pattern of node connections, this geometric construct can be viewed as a spacetime version of a diamond crystal. The ``covalent bonds'' of this ``crystal'' are lightray segments, i.e. causal links.

\section{Discretizations of curved spacetime} \label{sec:curved}

The work reported in this paper arose from a desire to use null 3d boundaries in discrete models of quantum gravity. In this section, we briefly discuss the utility of the geometric constructs discussed above for discretizations of curved spacetime. 

It is well-known that a (pseudo-)Riemannian space of any dimension can be approximated by a triangulation consisting of flat simplices. This is the basis of Regge calculus \cite{Regge:1961px,Williams:1991cd}. There, curvature is encoded in the deficit angles that arise when the simplices are glued together. Unfortunately, such a scheme \emph{cannot} work with tessellations of spacetime into null-faced 4-polytopes. First, as we've seen in section \ref{sec:polytopes:simplex}, one cannot triangulate \emph{any} spacetime with null-faced 4-simplices, due to the combinatorics of the faces' causal classification. One could still hope to tessellate spacetime with more complicated null-faced 4-polytopes. For flat spacetime, this is certainly possible, e.g. by tiling with parallelotopes. However, it turns out that null polyhedra cannot capture all the necessary kinds of curvature. The reason is the zero-sum condition \eqref{eq:s_diff} on the face areas. Because of this, a null hypersurface tessellated into flat polyhedra must have a constant area along all spacelike sections. This rules out even the simplest curved hypersurfaces in flat spacetime, e.g. lightcones. Since Ricci curvature is related to area increase via Raychaydhuri's equation, we conclude that it cannot be properly captured by a tessellation with null polyhedra.

One may get around this problem by using discretizations with \emph{curved} polyhedral cells. A direction that appears more promising is to generalize the lightray-threaded lattice from section \ref{sec:tiling}. In order to capture curvature, one may have to allow the cross-ratios of the lightrays' directions to vary, as well as the affine spacing of the nodes along them. Also, the number of ray segments connected to each node should be reduced. In particular, we want to allow lightrays to terminate, since they may fall through black hole horizons. A particularly promising construction is to leave five ray segments attached to each node - four past-going and one future-going, or vice versa. As before, the single future-going (say) segment is a continuation of one of the other four. The local variables will be the inner products of the four past-going (or the four future-going) segments at each node. The lone future-going (past-going) one is understood to have the same scalar products as its past-going (future-going) counterpart, but with opposite sign. This setup has a number of desirable features:
\begin{itemize}
 \item Each node is connected to five neighbors. The same is true for the dual of a 4d triangulation. This is encouraging, since triangulations are a successful discretization scheme, both in classical Regge calculus \cite{Gentle:2002ux} and in causal dynamical triangulation models \cite{Loll:2000my,Ambjorn:2010rx}.
 \item The 6 scalar products among the 4 independent segments at each node, together with the segments' null nature, are just enough to determine the metric at the node. This is consistent with the analogy between our nodes and dual triangulation sites. Indeed, in Regge calculus, the spacetime metric is defined within each 4-simplex, i.e. at each node of the dual triangulation.
 \item At each node, the causal pattern of ray segments is such that we have a null hyperplane acting as a future (past) causal boundary: all segments are either in it, or lie in its past (future). This is the hyperplane defined by the continuous lightray passing through the node, i.e. by the two paired past-going and future-going segments. Having a causal boundary at each node is useful for two purposes: to delineate a null outer boundary of the spacetime region of interest, and to allow for the dynamical appearance of horizons within that region.
\end{itemize}
We expect the future study of this structure to be rewarding.

\section*{Acknowledgements}		

I am grateful to Eugenio Bianchi for a discussion at Perimeter Institute. This work is supported in part by the NSF grant PHY-1205388 and the Eberly Research Funds of Penn State. The early stages were carried out in Tel Aviv University, supported in part by the Israeli Science Foundation center of excellence, the US-Israel Binational Science Foundation (BSF), the German-Israeli Foundation (GIF) and the Buchmann Scholarship Fund.

\end{document}